\documentclass[12pt]{article}

\usepackage[a4paper, textwidth=17cm, textheight=25.5cm, left=2cm]{geometry}

\usepackage{amsmath,amsthm,amssymb}
\allowdisplaybreaks[3]

\numberwithin{equation}{section}

\newtheorem{Definition}{Definition}
\newtheorem{theorem}{Theorem}[section]

\theoremstyle{remark}
\newtheorem{Remark}[theorem]{Remark}
\newcommand\CB{{\mathcal B}}
\newcommand\R{{\mathbb R}}
\newcommand\CR{{ \mathcal R}}

\newcommand\X{{\R^d}}
\newcommand\N{{\mathbb N}}
\newcommand\M{{\mathcal M}}

\newcommand\K{{\mathbb  K}}

\newcommand \Y{{\R^*_+ \times \X}}

\begin{document}

\title{Applied philosophy in mathematics}
\author{ Yu.G. Kondratiev\\
Bielefeld University, Germany\\
and Dragomanov University, Kyiv, Ukraine\\
}
%\date{}

%\begin{document}

\maketitle

{\small

\begin{center}
{\bf Abstract}

We show a possibility to apply certain philosophical concepts to
the analysis of concrete mathematical structures. Such application
gives a clear justification of topological and geometric properties
of considered mathematical objects. 

\end{center}
\noindent \vspace{2mm}

{\bf Keywords: discrete measures, configuration spaces, reflection maps} \vspace{2mm}

{\bf {MSC:  20C99, 2005, 47B38  }} 		 

%\newpage

\section{Introduction}

Interdisciplinary studies represent one of main trends in the modern science. Several essential problems in the science
and its applications need combination and interaction of methods and ideas from different areas of our knowledge.
But there appear many practical difficulties in the realization of an interdisciplinary approach. The point is that experts in
a particular topic may be not so deep involved in related areas outside their competence For example, we observe
very active development of mathematical modelings in the biology and ecology. But the use of such models by experts
from these disciplines is essentially restricted by the lack of mathematical techniques. From the other hand side,
mathematical models looks as very simplified and degenerated ones for experts in life sciences. The only way to
overcome these difficulties is to create the practical and patient collaboration  between concrete scientists.

Another traditional and old circle of discussions (and many speculations) concerns the relation between concrete 
sciences and the philosophy.  In the time of Newton and Leibniz the concept of the Naturphilosophy was
a commonly accepted basis for the unification of several scientific disciplines. But necessary specialization
and dissipation of particular sciences  produced the divergence of philosophy and concrete sciences and even
certain moral prejudices. No doubts, concrete results in physics, biology etc. are
still very stimulating for philosophical studies. But we would like to show that there 
exists a fruitful inverse influence. The aim of this work is to illustrate a natural applied
aspects of particular philosophical concepts in the framework of the mathematics. We did choose a concrete 
mathematical object for this illustration. Due to the interdisciplinary character of this journal, we restricted ourselves 
at  few basic observations about this object. Our explanations with necessity will be restricted technically to as less advanced level
as it is possible to keep an interest of  not only especially mathematical audience. For detailed mathematical description
of related structures we refer to \cite{FKKO 19},\cite{FKKO1 19}.

The myth of Plato's Cave served as one of the motivations for creating his concept about  the world of ideas and the world of things. In the dialogue ''State''  he gives a number of examples  illustrating this position. As we know, Plato considered mathematics as one of the most important sections, used in the construction of his philosophical system. Mathematical theories can serve as simple and illustrative tools for the existence of a "world of ideas" and a "world of things."  In a number of model situations, we are dealing with objects that appeal from our observations in physics, biology, ecology etc.  But full understanding of the mathematical structures of these models in many cases requires consideration of more general mathematical theories, which under some canonical mapping lead to the considered model situations. 

As an example, we can cite a number of recent works  on the study of spaces of random discrete measures. Such measures arise in many applications, in particular, in the theory of representations of current groups (Gelfand-Graev-Vershik), in models of biosphere (motivated by V. Vernadsky), etc. It turned out that the correct understanding of the topology and geometry of spaces of discrete measures naturally arises from the suitable configuration spaces on which these concepts are already well known.These configuration spaces we called Plato spaces and their elements are interpreted as "mathematical ideas" for our models of observed phenomena (of things). Maybe a naive illustration of this approach is related to Manin's concept of the adelic world as the space of Ideas \cite{Manin}. Moreover, the real component of adeles can be regarded as a "shadow" in the sense of Plato's theory. Number of similar examples in mathematical models can be big. Below we will describe a realization of the mentioned concept in a particular
case of random discrete measures.

Configuration spaces form an important and actively developing area in the infinite dimensional analysis.
From one hand side, these spaces represent reach mathematical structures which combine in a very non-trivial
way continuous and combinatoric aspects of the analysis. From the another side, configuration spaces
give natural  mathematical techniques in the applications to problems of mathematical physics, biology and ecology.

Spaces of discrete Radon measures  (DRM) may be considered as generalizations of configuration spaces. Main specific 
moment in the study of these spaces is such that the supports of discrete measures are typically not more
configurations. The latter change drastically technical methods in their study. Note that spaces of DRM have several
motivations coming from different areas of mathematics and applications, see comments below.

When choosing a model, one needs to take into account different features which are relevant for the behavior and properties of the system.
The considered state space can be chosen as a discrete set  or continuous, such as $\R^d$ or more generally, a Riemannian manifold $X$. While discrete models are easier to analyze (e.g. \cite{MR776231}) and yield more results, a continuous state space models a physical system more realistically.
  
  Bounded region vs. unbounded region/state space: A bounded region makes more sense from a modeling point of view. On the other hand, one needs to take into account the interaction of particles with the boundary. A way to circumvent this is by considering an unbounded region and restricting the system after analyzing the model. The kind of region also determines whether a finite or an infinite amount of particles should be considered.
   Another advantage of an unbounded region with an infinite number of particles is that phase transitions may be observed since invariant measures may not be uniquely determined. For examples, see \cite{MR3421803} and the references therein.

  Different mechanisms yield different behaviors of the system. This choice of course depends on the desired phenomenon which is to be modeled. 
There are some additional options which were already mentioned above. For our situation, we choose a specific version of a continuous particle system with unbounded state space $\R^d$. Instead of considering a homogeneous configuration space, the particle system comes from the cone of positive discrete Radon measures. One specific property of this object is that particles in the space $\R^d$ are assigned a positive number, or "mark'', which represents a property of the particle such as weight. Some general analytic and geometric considerations for models on the cone of Radon measures have been carried out in \cite{MR3708379,MR3394626}.

Note that this approach differs from the so-called marked configuration spaces considered in \cite{2006math......8347K,2006math......8344A}. On the other hand, there is a direct relation to the extended configuration space  which we describe
below. While the analysis and dynamics on the cone are of special interest and the modeling possibilities of the cone are useful in applications, one may also give some motivations for this object without referring to these analytical properties or configuration spaces in general. Below we explains three motivations from theoretical biology, probability theory and representation theory.

The mathematical object of interest for us is the cone of positive DRM, defined by
\begin{displaymath}
 \K(\X):=\left\{\eta=\sum_i s_i\delta_{x_i}\in\M(\X)\middle| s_i\in (0,\infty), x_i\in\X\right\}
\end{displaymath}
where by convention, the zero measure $0\in\K(\X)$ is included. This work is concerned with the analytic properties of the cone. On the other hand, there are three approaches which justify the use of this object without even considering its analytical properties. For one, there is the aspect of modeling biological systems. Second, the cone appears naturally when considering certain generalized stochastic processes. Third, the cone is given as the space where Gamma measures are localized, which emerge from representation theory for current groups. These three motivations will be explained  below.

There is an external non-mathematical motivation to study particle systems realized as elements of the cone. Namely, Vladimir Vernadsky  wrote the following:
\begin{itemize}
 \item ``Organisms [...] are always separated from the surrounding inert matter by a clear and firm boundary.'' \cite[p. 56]{Vernadsky1998}
 \item ``Living matter [...] is spread over the entire surface of the Earth in a manner analogous to a gas [...].''\cite[p. 59]{Vernadsky1998}
 \item ``In the course of time, living matter clothes the whole terrestrial globe with a continuous envelope [...].'' \cite[p. 60]{Vernadsky1998}
\end{itemize}
This can be interpreted in the sense that system of living matter should possess two properties: For one, the system should have a discrete nature. Furthermore, there is living matter everywhere in the system. In mathematical terms, this means that the support of this system should be dense in the underlying position space. Lastly, to be realistic, the system should have finite local mass due to the physical limitations of our world. The mathematical realization of these properties is given by the cone.

The second motivation comes from the theory of generalized stochastic processes, i.e. processes on the space $\mathcal{D}^\prime(\X)$ of generalized functions. By \cite[Thm. 3.3.24]{MR1155400}, infinitely divisible processes on $\mathcal{D}^\prime(\X)$ are actually concentrated on the subspace $\K(\X)$. Note that this result holds independently of the topological and analytical considerations done in later chapters. For a subclass of measures, the so-called Gamma measures, we will also show a direct proof of this statement.

Measures supported on $\K(\X)$ naturally appear in the study of representations for current groups. Namely, when studying so-called commutative models of representations of $(SL(2,\R))^{\X}$. When considering representations with respect to the unipotent subgroup of $(SL(2,\R))^{\X}$, we arrive at spectral measures which are defined on the space $\mathcal{D}^\prime(\X)$ and supported on $\K(\X)$. Furthermore, these measures show some invariance properties. These considerations were first done by Gelfand, Graev and Vershik \cite{MR829048}. Later, Tsilevich, Vershik and Yor \cite{MR1853759} used this as a starting point to further analyze so-called Gamma processes.

As seen here, these measures supported on the cone $\K(\X)$ appear naturally without any \emph{a priori} restriction of the spaces or aspects of modeling.

There is another mathematical explanation why it makes sense to consider $\K(\X)$. If we take the class of Gamma-Poisson-measures on the extended configuration space $\Gamma(\R^*_+\times\X)$, we see that these measures assign full mass to the subset of configurations with finite local mass, or Plato configurations. These configurations can be identified with objects in the cone, i.e. there exists a one-to-one correspondence between the so-called Plato space $\Pi(\R^*_+\times\X)$ and the cone $\K(\X)$
\cite{FKKO  19}.

\section{Preliminaries}\label{prelims}

This section will include basic concepts  from configuration spaces $\Gamma(\Y)$, the cone $\K(\X)$ and the connection between these two.

\subsection{The Cone of Positive Discrete Radon Measures}\label{the_cone}

We start by the introduction of the cone of positive discrete Radon measures
as the subset of the space of Radon measures  $\M(\X)$. Furthermore, the notion of the support of a measure and relations between elements in $\K(\X)$ are defined. Recall that by Vernadsky's theory of living matter, a system should be dense everywhere, discrete and have finite local mass.

One more property which we want from our system is that its elements are indistinguishable in the sense that the system given by $(s_i,x_i)_{i\in I}$ and $(s_{\pi(i)},x_{\pi(i)})_{i\in I}$ behave the same, where $I$ is some countable index set and $\pi$ an arbitrary permutation of $I$. One possibility is to realize our system as sums of point masses $\delta_y$, where $y$ is either the mark and position, or just the position of a particle, depending on the setup. This automatically yields a discrete particle system. To obtain the other two properties, it is useful to let $y$ represent the position of a particle, while the mark is considered as a weight of the point mass. These properties become clear when we consider a certain class of measures, namely, Gamma measures. 

\begin{Definition}\label{def22}
 \begin{enumerate}
  \item The cone of nonnegative discrete Radon measures is defined as follows:
  \begin{displaymath}
   \K(\X):=\left\{\eta=\sum_i s_i\delta_{x_i}\in\M(\X)\middle| s_i\in (0,\infty), x_i\in\X\right\}
  \end{displaymath}
  By convention, the zero measure $\eta=0$ is included in $\K(\X)$.
  \item We denote the support of $\eta\in\K(\X)$ by
  \begin{displaymath}
   \tau(\eta):=\{x\in\X\mid 0<\eta(\{x\})=:s_x(\eta)\}.
  \end{displaymath}
  If $\eta$ is fixed, we write $s_x:=s_x(\eta)$.
  \item For $\eta,\xi\in\K(\X)$ we write $\xi\subset\eta$ if $\tau(\xi)\subset\tau(\eta)$ and $s_x(\xi)=s_x(\eta)$ for all $x\in\tau(\xi)$. If additionally $|\tau(\xi)|<\infty$, we write $\xi\Subset\eta$.
  \item For a function $f\in C_c(\X)$, denote the pairing with an element $\eta\in\K(\X)$ by
  \begin{displaymath}
   \langle f,\eta\rangle:=\sum_{x\in\tau(\eta)}s_x f(x).
  \end{displaymath}
 \end{enumerate}
\end{Definition}

While $\K(\X)$ can be viewed as a subset of the space of positive Radon measures $\M(\X)$, it is not advisable to consider it as a subset topologically. This method works for the space $\Gamma(Y)$ introduced below, as will be explained later. For $\K(\X)$, it does not yield satisfactory topological results. Instead, we keep Plato's theory in mind and see $\K(\X)$ as the real-world projection of another space, called the Plato space $\Pi(\Y)$. 

\subsection{Plato's theory}

As stated in the introduction, the cone $\K(\X)$ is a suitable object to describe particle systems in the real world. On the other hand, the question arises how to define and interpret mathematical structures on the space $\K(\X)$. As a motivation, we give a short overview of Plato's theory of forms.

In the theory, Plato stated that observations in the real world are mere projections of higher forms or ideas. One way to picture this is the so-called cave allegory, which was recited by Ross (1951) as follows: ``A company of men is imprisoned in an underground cave, with their heads fixed so that they can look only at the back wall of the cave. Behind them across the cave runs a wall behind which men pass, carrying all manner of vessels and statues which overtop the wall. Behind these again is a fire. The prisoners can only see the shadows [...] of the things carried behind the wall, and must take these to be the only realities'' \cite[P. 69]{ross_plato}.

Applied to our setting, the space $\K(\X)$ is interpreted as the shadows projected onto the cave wall. On the other hand, the space $\Pi(\Y)$ which will be introduced below is the space of forms or ideas, represented by the objects carried in front of the fire. While the space $\K(\X)$ is taken to be our reality, we use the space $\Pi(\Y)$ to define mathematical operations. The spaces are connected via the bijection $\CR\colon\Pi(\Y)\to\K(\X)$ introduced below. In accordance with the cave allegory, $\CR$ is also called reflection mapping.

\subsection{Configuration Spaces}\label{conf}

As we will see in the next chapter, the Plato space $\Pi(\Y)$ is a  specific subset of the so-called configuration space $\Gamma(\Y)$, which will fulfill the assumptions stated heuristically in Chapter \ref{the_cone}.

In general, the space of locally finite configurations $\Gamma(Y)$ is the space of all subsets of $Y$ which are finite in any compact set $\Lambda\subset Y$. The following definition makes this notion more precise.
\begin{Definition}
 Let $Y$ be a locally compact Hausdorff space. The space of locally finite configurations over $Y$ is defined as
 \begin{displaymath}
  \Gamma(Y)=\{\gamma\subset Y\colon |\gamma\cap\Lambda|<\infty\ \forall\Lambda\subset Y\ \mathrm{compact}\}
 \end{displaymath}
 where $|\cdot|$ denotes the number of elements of a set.
\end{Definition}
From a physical perspective, $Y$ is considered as phase space of an interacting particle system. A configuration $\gamma\in\Gamma(Y)$ represents a set of indistinguishable agents (e.g. particles, plants) which may interact with each other. In our considerations, we always consider $Y=\Y$. More generally, $\X$ could be replaced by some more general locally comapct space $X$. In this chapter, we recall some properties of $\Gamma(Y)$ which will form the basis for the Plato space $\Pi(\Y)$.

\subsubsection{Topology and Measurable Structure of $\Gamma(Y)$}\label{topology_gamma}

There exists a natural embedding of $\Gamma(Y)$ into the space of Radon measures $\M(Y)$ on $Y$, namely
\begin{displaymath}
 \Gamma(Y)\ni\gamma\mapsto\sum_{y\in\gamma}\delta_y\in\M(Y)
\end{displaymath}
where $\delta_y$ denotes the Dirac measure at point $y\in Y$.  We equip $\Gamma(Y)$ with the vague topology induced by $\M(Y)$, i.e. the coarsest topology such that the following mappings are continuous for all $f\in C_c(Y)$, where $C_c(Y)$ denotes the space of continuous functions with compact support:
\begin{displaymath}
 \Gamma(Y)\ni\gamma\mapsto\langle f,\gamma\rangle=\sum_{y\in\gamma}f(y)
\end{displaymath}
In fact, $\Gamma(Y)$ equipped with this topology is a Polish space. A more detailed analysis of the topological properties of $\Gamma(Y)$ can be found in \cite{MR2226411}.

The construction of a topology enables us to consider the Borel-$\sigma$-algebra $\CB(\Gamma(Y))$. It should be noted that this $\sigma$-algebra coincides with the $\sigma$-algebra generated by the following mappings:
\begin{displaymath}
 N_\Lambda\colon\Gamma(Y)\to\N_0, \gamma\mapsto N_\Lambda(\gamma)=|\gamma\cap\Lambda|,\ \Lambda\in\CB_c(Y)
\end{displaymath}
where $\CB_c(Y)$ denotes all pre-compact Borel subsets of $Y$, see e.g. \cite{MR1914839}.

We give another construction of the measurable space $(\Gamma(Y),\CB(\Gamma(Y))$ which will be useful for other considerations. For $\Lambda\in\CB_c(Y)$, we define the space of configurations supported in $\Lambda$.
\begin{displaymath}
 \Gamma(\Lambda):=\{\gamma\in\Gamma(Y)\colon\gamma\cap\Lambda=\gamma\}.
\end{displaymath}
 Furthermore, for $n\in\N$, consider the set of $n$-point-configurations supported in $\Lambda$:
\begin{displaymath}
 \Gamma^{(n)}(\Lambda):=\{\gamma\in\Gamma(\Lambda)\colon|\gamma|=n\}, \Gamma^{(0)}(\Lambda):=\{\emptyset\}
\end{displaymath}
Since $\gamma\in\Gamma(Y)$ is locally finite, the elements of $\Gamma(\Lambda)$ are finite and we have the disjoint decomposition
\begin{equation}\label{decomp}
 \Gamma(\Lambda)=\bigcup_{n=0}^\infty\Gamma^{(n)}(\Lambda).
\end{equation}
We can represent $\Gamma^{(n)}(\Lambda)$ via symmetrization of the underlying space:
\begin{equation}\label{sym}
 \tilde{\Lambda}^n/S_n\simeq\Gamma^{(n)}(\Lambda)
\end{equation}
where
\begin{displaymath}
 \tilde\Lambda^n:=\{(x_1,\dotsc,x_n)\in\Lambda^n\mid x_i\neq x_j\ \forall i\neq j\}
\end{displaymath}
the off-diagonals and $S_n$ the symmetric group of $n$ elements.
% Denote the bijection $ \eqref{sym}$  by $\sym_n$. 
 This way, $\Gamma^{(n)}(\Lambda)$ can be equipped with the topology induced via $\Lambda^n$. Furthermore, $\Gamma(\Lambda)$ is equipped with the topology of disjoint unions. Hence, we can define the Borel-$\sigma$-algebra $\CB(\Gamma(\Lambda))$ given by this topology.

For two sets $\Lambda_1,\Lambda_2\in\CB(Y), \Lambda_2\subset\Lambda_1$, define the projection mapping
\begin{displaymath}
 p_{\Lambda_1,\Lambda_2}\colon\Gamma(\Lambda_1)\to\Gamma(\Lambda_2), \gamma\mapsto\gamma\cap\Lambda_2
\end{displaymath}
where we set $p_{\Lambda_2}:=p_{Y,\Lambda_2}$. It was shown in e.g. \cite{MR1295944} that $(\Gamma(Y),\CB(\Gamma(Y))$ is the projective limit of the spaces $(\Gamma(\Lambda),\CB(\Gamma(\Lambda))$ for $\Lambda\in\CB_c(Y)$. This especially implies that the mappings $p_\Lambda$ are $\CB(\Gamma(Y))$-$\CB(\Gamma(\Lambda))$-measurable. The construction of $\CB(\Gamma(Y))$ via projections will play an important role in the construction of measures on $\Gamma(Y)$.

\subsubsection{The Space of Finite Configurations}\label{sec_finite_conf}

For mathematical purposes, it is important to also consider the space $\Gamma_0(Y)$ of finite configurations, i.e.
\begin{displaymath}
 \Gamma_0(Y):=\{\gamma\in\Gamma(Y)\colon|\gamma|<\infty\}
\end{displaymath}
where $|\cdot|$ denotes the number of elements of a set. While the definition implies that $\Gamma_0(Y)$ is a subset of $\Gamma(Y)$, the interpretation is a different one: $\Gamma_0(Y)$ serves as a mathematical counterpart to the physical space $\Gamma(Y)$. Also, the spaces $\Gamma(Y)$ and $\Gamma_0(Y)$ are topologically different: While $\Gamma(Y)$ is seen as a subspace of $\M(Y)$ with the inherited topology, we use a different approach for $\Gamma_0(Y)$ which will be explained in this chapter. The approach is similar to the one used in Chapter \ref{topology_gamma}, but yields different results. We set
\begin{displaymath}
 \Gamma_0^{(n)}(\Lambda):=\Gamma^{(n)}(\Lambda)
\end{displaymath}
where $\Lambda$ is an arbitrary Borel subset of $Y$. Since we only deal with finite configurations, we may use decomposition \eqref{decomp} for $\Lambda=Y$, i.e.
\begin{displaymath}
 \Gamma_0(Y)=\bigsqcup_{n=0}^\infty\Gamma^{(n)}_0(Y).
\end{displaymath}
Furthermore, we may consider the symmetrization \eqref{sym} to obtain
\begin{displaymath}
 \tilde{Y}^n/S_n\simeq\Gamma^{(n)}(Y).
\end{displaymath}
For $\Gamma^{(n)}(Y)$, we choose the topology induced by the space $Y^n$. For $\Gamma_0(Y)$, we may use the topology of disjoint unions. For a more detailed description of the topology used here, we refer to \cite{MR1914839}.

\begin{Remark}
 The purpose of the space of finite configurations will become clearer once we examine specific models. Since the models are introduced on the cone, we postpone this discussion until after we have introduced the relevant spaces related to $\K(\X)$.
\end{Remark}

\subsection{Relation Between $\K(\X)$ and $\Gamma(\Y)$: The Plato Space $\Pi(\Y)$}\label{relation_pi_k}

In this section, we want to establish the connection between the configuration space $\Gamma(\R^*_+\times\X)$ and the cone $\K(\X)$. Our goal is to define a certain subspace $\Pi(\Y)\subset\Gamma(\Y)$ such that there exists a one-to-one-correspondence between $\Pi(\Y)$ and $\K(\X)$ in the following form:
\begin{displaymath}
 \CR\colon\Pi(\Y)\to\K(\X),\gamma=\sum_{(s,x)\in\gamma}\delta_{(s,x)}\mapsto\sum_{(s,x)\in\gamma}s\delta_x.
\end{displaymath}
In terms of Plato's theory, this mapping takes ideas $\gamma\in\Pi(\Y)$ and projects (or reflects) them to real-world objects $\eta\in\K(\X)$. Obviously, $\CR$ is not defined on the whole space $\Gamma(\R^*_+\times\X)$. Therefore, we need to construct a suitable subspace. In other terms, the Plato space constructed below is also known as the set of pinpointing configurations with finite local mass, denoted by $\Gamma_\mathrm{pf}(\R^*_+\times\X)$. We explore these two properties in more detail below.

Define the set of pinpointing configurations $\Gamma_\mathrm{p}(\Y)\subset\Gamma(\Y)$ as all configurations such that if $(s_1,x_1),(s_2,x_2)\in\gamma$ with $x_1=x_2$, then $s_1=s_2$.
\begin{Remark}
  The pinpointing property ensures that there are no two elements of a system at the same position. Due to the shape of elements in $\K(\X)$, it is obvious that this would not be possible.
\end{Remark}
Let us now take into account the second property of $\Pi(\Y)$. To this end, we define the local mass of a configuration.
\begin{Definition}
 For a configuration $\gamma\in\Gamma_\mathrm{p}(\R^*_+\times\X)$ and $\Lambda\subset\X$ compact, set the local mass as
 \begin{displaymath}
  \gamma(\Lambda)=\int_{\R_+ \times \X} s {1}_{\Lambda} (x)   d\gamma(s,x)=\sum_{(s,x)\in\gamma}s {1}_ {\Lambda}(x)\in[0,\infty]
 \end{displaymath}
\end{Definition}
This notion enables us to define the Plato space.
\begin{Definition}
 The Plato space $\Pi(\Y)\subset\Gamma(\Y)$ is defined as the space of all pinpointing configurations with finite local mass, i.e.
 \begin{displaymath}
  \Pi(\Y):=\Gamma_\mathrm{pf}(\Y)=\{\gamma\in\Gamma_p\mid\gamma(\Lambda)<\infty\text{ for all }\Lambda\subset\R^d\text{ compact}\}.
 \end{displaymath}
\end{Definition}
\begin{Remark}
 \begin{enumerate}
  \item The property of finite local mass accounts for the third property stated in Chapter \ref{the_cone}. It ensures that the system only has finite mass in any bounded volume, which makes it physically viable.
  \item The pinpointing property as well as the finiteness of local mass are sufficient to make $\CR\colon\Pi(\Y)\to\K(\X)$ bijective.
  \item The state space needs to be of the specific form $Y=\R^*_+\times X$ for the notion of pinpointing configurations to make sense.
 \end{enumerate}
\end{Remark}

\begin{Definition}
 Let $f\in C_c(\Y)$ and $\eta\in\K(\X)$. Define the following pairing:
 \begin{displaymath}
  \langle\langle f,\eta\rangle\rangle:=\langle f,\CR^{-1}\eta\rangle=\sum_{(s,x)\in\CR^{-1}\eta}f(s,x)
 \end{displaymath}

\end{Definition}

\section{Topology and Measure-Theoretical Structures}\label{sec:topology}

In this chapter, we want to introduce a suitable topology on the cone $\K(\X)$. To this end, we consider the topology induced on $\Pi(\Y)$ by the extended configuration space $\Gamma(\Y)$. Next, we use the mapping $\CR$ to induce a topology on $\K(\X)$.

\subsection{Topology on the Cone $\K(\X)$.}\label{topol_k}

The Plato space $\Pi(\Y)$ naturally inherits the topological structure of $\Gamma(\Y)$, i.e. the topology is given by the vague topology induced from the space of Radon measures $\M(\Y)$. For a detailed description of topological and metric characterizations, see e.g. \cite{MR2226411}.

\begin{Remark}
 The space $\Pi(\Y)$ is not complete: Take for example some $x_0\in\X$ and $s_1\neq s_2\in\R^*_+$. Furthermore, consider sequences $s_i^{(n)},x_i^{(n)}, i=1,2$ with $s_1^{(n)}\neq s_2^{(n)}, x_1^{(n)}\neq x_2^{(n)}$ for all $n\in\N$ and
\begin{displaymath}
 s_i^{(n)}\to s_i, x_i^{(n)}\to x_i,\ n\to\infty, i=1,2.
\end{displaymath}
Set
\begin{align*}
 \gamma^{(n)}&:=\{(s_1^{(n)},x_1^{(n)}),(s_2^{(n)},x_2^{(n)})\}\in\Pi(\Y)
 \\
 \gamma&:=\{(s_1,x_0),(s_2,x_0)\}\in\Gamma(\Y)\setminus\Pi(\Y)
\end{align*}
Let $f\in C_c(\Y)$. Then
\begin{align*}
 |\langle f,\gamma^{(n)}\rangle-\langle f,\gamma\rangle|&=|f(s_1^{(n)},x_1^{(n)})+f(s_2^{(n)},x_2^{(n)})-f(s_1,x_0)-f(s_2,x_0)|
 \\
 &\leq|f(s_1^{(n)},x_1^{(n)})-f(s_1,x_0)|+|f(s_2^{(n)},x_2^{(n)})-f(s_2,x_0)|
 \\
 &\to 0,\ n\to\infty.
\end{align*}
Therefore, $\gamma^{(n)}\to\gamma,\ n\to\infty$ in $\Gamma(\Y)$ and $\Pi(\Y)$ is not complete.
\end{Remark}

From a naive point of view, it seems to make sense to consider the embedding $\K(\X)\subset\M(\X)$ of the cone into the space of Radon measures, equipped with the vague topology. Unfortunately, this topology has no relation to the vague topology introduced above on $\Pi(\Y)$. In the spirit of Plato's theory of ideas, the connection between $\Pi(\Y)$ and $\K(\X)$ is essential. Therefore, we consider the final topology on $\K(\X)$ induced by the reflection mapping $\CR$, i.e. the finest topology such that the mapping
\begin{displaymath}
 \CR\colon\Pi(\Y)\to\K(\X), \gamma=\sum_{(s_x,x)\in\gamma}\delta_{(s_x,x)}\mapsto\sum_{x\in\tau(\gamma)}s_x\delta_x
\end{displaymath}
is continuous. Here, we set for $\gamma\in\Pi(\Y)$,
\begin{displaymath}
 \tau(\gamma):=\{x\in\X\mid\exists s\in\R^*_+\colon(s,x)\in\gamma\}
\end{displaymath}
the support of $\gamma$. The usage of this topology has the obvious side effect that $\CR$ becomes a homeomorphism, which is helpful in and of itself in other regards. 
%\subsection{Measures on $\Pi(\Y)$ and $\K(\X)$}\label{meas_pi_k}

In the further development of the considered theory is important to implement the construction of a  class of probability measures on $\Pi(\Y)$, namely, Poisson measures. The construction may be done on the larger space $\Gamma(\Y)$. For the class of Poisson measures, we  can show that they assign full mass to the Plato space $\Pi(\Y)$.
To obtain measures on $\K(\X)$, we use the pushforward of measures on $\Pi(\Y)$ under the mapping $\CR$. A certain subclass of specific interest is the class of Gamma measures. For the detailed analysis we refer the reader to \cite{FKKO 19}, \cite{FKKO1 19}.


\begin{thebibliography}{10}

\bibitem{2006math......8344A}
S.~{Albeverio}, Y.~G. {Kondratiev}, E.W. {Lytvynov}, and G.~F. {Us}.
\newblock {Analysis and geometry on marked configuration spaces}.
\newblock {\em arXiv Mathematics e-prints}, page math/0608344, Aug 2006.


\bibitem{MR1612725}
S.~Albeverio, Yu.~G. Kondratiev, and M.~R\"{o}ckner.
\newblock Analysis and geometry on configuration spaces.
\newblock {\em J. Funct. Anal.}, 154(2):444--500, 1998.

\bibitem{MR3421803}
A.~Daletskii, Y.~G. Kondratiev, and Y.~Kozitsky.
\newblock Phase transitions in continuum ferromagnets with unbounded spins.
\newblock {\em J. Math. Phys.}, 56(11):113502, 16, 2015.

\bibitem{FKKO 19} D.Finkelshtein, Y.Kondratiev, P.Kuchling, and
M.J.Olivera.
\newblock Analysis and geometry on the cone of discrete
Radon measures I
\newblock {\em Methods Funct. Anal. Topology}. to appear, 2019

\bibitem{FKKO1 19} D.Finkelshtein, Y.Kondratiev, P.Kuchling, and
M.J.Olivera.
\newblock Analysis and geometry on the cone of discrete
Radon measures II
\newblock {\em Methods Funct. Anal. Topology}. to appear, 2019


\bibitem{MR829048}
I.~M. Gelfand, M.~I. Graev, and A.~M. Vershik.
\newblock Models of representations of current groups.
\newblock In {\em Representations of {L}ie groups and {L}ie algebras
  ({B}udapest, 1971)}, pages 121--179. Akad. Kiad\'{o}, Budapest, 1985.

\bibitem{MR3708379}
D.~Hagedorn, Y.~G. Kondratiev, E.~Lytvynov, and A.~Vershik.
\newblock Laplace operators in gamma analysis.
\newblock In {\em Stochastic and infinite dimensional analysis}, Trends Math,
  pages 119--147. Birkh\"{a}user/Springer, [Cham], 2016.

\bibitem{MR3041709}
D.~Hagedorn, Y.~G. Kondratiev, T.~Pasurek, and M.~R\"{o}ckner.
\newblock Gibbs states over the cone of discrete measures.
\newblock {\em J. Funct. Anal.}, 264(11):2550--2583, 2013.

\bibitem{MR1914839}
Y.~G. Kondratiev and T.~Kuna.
\newblock Harmonic analysis on configuration space. {I}. {G}eneral theory.
\newblock {\em Infin. Dimens. Anal. Quantum Probab. Relat. Top.},
  5(2):201--233, 2002.

\bibitem{MR2226411}
Y.~G. Kondratiev and O.~Kutoviy.
\newblock On the metrical properties of the configuration space.
\newblock {\em Math. Nachr.}, 279(7):774--783, 2006.


\bibitem{MR3394626}
Y.~G. Kondratiev, E.~Lytvynov, and A.~Vershik.
\newblock Laplace operators on the cone of {R}adon measures.
\newblock {\em J. Funct. Anal.}, 269(9):2947--2976, 2015.

\bibitem{2006math......8347K}
Y.~G. {Kondratiev}, E.~W. {Lytvynov}, and G.~F. {Us}.
\newblock {Analysis and geometry on $R_+$-marked configuration spaces}.
\newblock {\em arXiv Mathematics e-prints}, page math/0608347, Aug 2006.

\bibitem{MR776231} 
T.M.~Liggett,
\newblock {\em Interacting particle systems.}
\newblock {Grundlehren der Mathematischen Wissenschaften [Fundamental Principles of Mathematical Sciences]}
\newblock {276. Springer-Verlag, New York, 1985. {\rm xv}+488 pp.}

\bibitem{Manin} 
Yu.I.~Manin
\newblock {Reflections on arithmetic physics}
\newblock{ \em Invariance and string theory}
\newblock Academic Press, 1989, 293-303


\bibitem{ross_plato}
W.~D. Ross.
\newblock {\em Plato's theory of ideas}.
\newblock Clarendon Pr., Oxford, 1951.

\bibitem{MR1295944}
H.~Shimomura. 
\newblock  Poisson measures on the configuration space and unitary representations of the group of diffeomorphisms.
\newblock {\em J. Math. Kyoto Univ.}, 34(3), 599--614, 1994. 

\bibitem{MR1155400}
A~Skorohod.
\newblock {\em Random processes with independent increments}, volume~47 of {\em
  Mathematics and its Applications (Soviet Series)}.
\newblock Kluwer Academic Publishers Group, Dordrecht, 1991.
\newblock Translated from the second Russian edition by P. V. Malyshev.

\bibitem{MR1853759}
N.~Tsilevich, A.~Vershik, and M.~Yor.
\newblock An infinite-dimensional analogue of the {L}ebesgue measure and
  distinguished properties of the gamma process.
\newblock {\em J. Funct. Anal.}, 185(1):274--296, 2001.

\bibitem{Vernadsky1998}
V.~I. Vernadsky.
\newblock {\em Living Matter in the Biosphere}, pages 56--60.
\newblock Springer New York, New York, NY, 1998.

\end{thebibliography}
\end{document}